\documentclass[12pt,twoside]{amsart}

\usepackage{amsfonts,amssymb}
\usepackage{xr}
\usepackage{graphicx}
\usepackage{amsmath}
\usepackage{epsfig}
\usepackage{color}
\def\R{\mathbb{R}}
\def\N{\mathbb{N}}

\def\epsilon{\varepsilon}

\def\tilde{\widetilde}

\def\trait (#1) (#2) (#3){\vrule width #1pt height #2pt depth #3pt}
\def\fin{\hfill\trait (0.1) (5) (0) \trait (5) (0.1) (0) \kern-5pt 
\trait (5) (5) (-4.9) \trait (0.1) (5) (0)}

\newcommand{\be}{\begin{equation}}
\newcommand{\ee}{\end{equation}}
\newcommand{\baa}{\begin{array}}
\newcommand{\eaa}{\end{array}}
\newcommand{\ba}{\begin{eqnarray}}
\newcommand{\ea}{\end{eqnarray}}

\newtheorem{theo}{\bf Theorem}[section]
\newtheorem{lem}[theo]{\bf Lemma}

\title[precised hardy inequalities]{ A note on precised Hardy inequalities on Carnot groups and Riemannian manifolds }
\author{Emmanuel Russ}
\author{Yannick Sire}

\begin{document}

\maketitle

\begin{abstract} 
We prove non local Hardy inequalities on Carnot groups and Riemannian manifolds, relying on integral representations of fractional Sobolev norms. 
\end{abstract} 

\noindent{\small{{\bf AMS numbers 2000: }}} Primary: 46E35. Secondary:  35R11, 42B35, 58J35.

\noindent{\small{{\bf Keywords: }}} Hardy inequalities, Lie groups, Riemannian manifolds, fractional powers.
\tableofcontents

%%%%%%%%%%%%%%%%%%%%%%%%%%%%%%%%%%%%%%%%
%%%%%%%%%%%%%%%%%%%%%%%%%%%%%%%%%%%%%%%%
\section{Introduction}
 In the whole paper, when two quantities $A(f)$ and $B(f)$ depend on a function $f$ ranging in some space $V$, the notation 
$$A(f)\lesssim B(f)\ \forall f\in V
$$ 
means that there exists $C>0$ such that $A(f)\leq CB(f)$ for all $f\in V$. Moreover, the notation
$$
A(f)\sim B(f)\ \forall f\in V
$$
means that there exists $C>0$ such that $C^{-1}B(f)\leq A(f)\leq CB(f)$ for all $f\in V$.
\subsection{Setting of the problem}

The simplest Hardy inequality on $\R^n$ asserts that, if $n\geq 3$, 
\begin{equation} \label{hardyrn}
\int_{\R^n} \frac{u^2(x)}{\left\vert x\right\vert^2} dx\lesssim \int_{\R^n} \left\vert \nabla u(x)\right\vert^2dx=\left\Vert u\right\Vert_{\dot{H}^1(\R^n)},\ \forall u\in {\mathcal D}(\R^n).
\end{equation}
A non local version of \eqref{hardyrn} can be given, where the $\dot{H}^1$ norm in the right hand side is replaced by a $\dot{H}^s$ norm for $0<s<\frac n2$ (see \cite{BCG}):
\begin{equation} \label{hardyrnfrac}
\int_{\R^n} \frac{u^2(x)}{\left\vert x\right\vert^{2s}} dx\lesssim \left\Vert u\right\Vert_{\dot{H}^s(\R^n)}^2,\ \forall u\in {\mathcal D}(\R^n).
\end{equation}
When $0<s<1$, it is well-known (see for instance \cite{adams}) that  $\left\Vert u\right\Vert_{\dot{H}^s(\R^n)}$ can be represented by means of an integral quantity involving first order differences of $u$, and (\ref{hardyrnfrac}) can therefore be rewritten as
\begin{equation} \label{hardyrnfracbis}
\int_{\R^n} \frac{u^2(x)}{\left\vert x\right\vert^{2s}} dx\lesssim \iint_{\R^n\times \R^n} \frac{\left\vert u(x)-u(y)\right\vert^2}{\left\vert x-y\right\vert^{n+2s}}dxdy.
\end{equation}
These Hardy inequalities ({\it i.e.} the local and the non local version) were transposed to the framework of the Heisenberg group in \cite{BCG,BCX}. More precisely, in the Heisenberg group ${\mathcal H}^d$ ($d\geq 1$), the following Hardy inequality was established in \cite{BCX}:
\begin{equation} \label{hardyhd}
\int_{{\mathcal H}^d}  \frac{u^2(x)}{\rho^2(x)}dx\lesssim \left\Vert \nabla_{\mathcal H}u\right\Vert_2^2,\ \forall u\in {\mathcal D}({\mathcal H}^d),
\end{equation}
where $\rho(x)$ denotes the distance of $x$ to the origin and $\nabla_{\mathcal H}$ stands for the gradient associated to the vector fields $Z_1,...,Z_{2d}$ (see \cite{BCX} and the notations therein). The non local version of (\ref{hardyhd}), which was proven in \cite{BCG} (where it was derived from precised inequalities involving Besov norms) says that, for $0<s<d+1$,
\begin{equation} \label{hardyhdfrac}
\int_{{\mathcal H}^d}  \frac{u^2(x)}{\rho^{2s}(x)}dx\lesssim \left\Vert u\right\Vert_{\dot{H}^s}^2,\ \forall u\in {\mathcal D}({\mathcal H}^d).
\end{equation}
When $0<s<1$, an integral representation for the fractional Sobolev homogeneous norm was proven in \cite{CRT} (note that an analogous representation holds in any connected Lie group wth polynomial volume growth, and even in any unimodular Lie group if one works with the inhomogeneous version of this norm), so that \eqref{hardyhdfrac} can be rewritten as
\begin{equation} \label{hardyhdfracbis}
\int_{{\mathcal H}^d}  \frac{u^2(x)}{\rho^{2s}(x)}dx\lesssim \iint_{{\mathcal H}^d\times {\mathcal H}^d} \frac{\left\vert u(x)-u(y)\right\vert^2}{\rho(y^{-1}x)^{2d+2+2s}}dxdy.
\end{equation}

\par

\medskip

\noindent Hardy inequalities in {\it local} versions on more general Lie groups, namely Carnot groups, were obtained in \cite{kombe}. In the present paper, we establish non local versions of these Hardy inequalities on Carnot groups, in the spirit of \eqref{hardyrnfracbis} and \eqref{hardyhdfracbis}. We also investigate the similar problem on Riemannian manifolds.
 
%This paper is devoted to the derivation of non local Hardy inequalities on Lie groups and Riemannian manifolds. By ``non local'', we mean that the gradient in the right-hand side is replaced by a non local expression, which is reminiscent of a Gagliardo-type norm. Precised Hardy inequalities have been obtained in the case of the Heisenberg group in \cite{BCG}. In the present work, we give a different version involving a non local term in the right hand-side of the inequality. We also consider two other situations

%\begin{itemize}
%\item The case of general Carnot groups, for which we provide a new precised Hardy inequality. 
%\item The case of general Lie groups with polynomial growth. 
%\end{itemize}  

%Notice that the Carnot groups, since they are stratified, enjoy polynomial growth. 

\subsection{The case of Lie groups}
We now describe the general framework for Lie groups. Let $G$ be a unimodular connected Lie group endowed with the Haar measure. By ``unimodular'', we mean that the Haar measure is left and right-invariant.  If we denote by $\mathcal  G$  the Lie algebra  of $G$, we consider a family 
$$\mathbb X= \left \{ X_1,...,X_k \right \}$$
of left-invariant vector fields on $G$ satisfying the H\"ormander condition, i.e. $\mathcal G$ is the lie algebra generated by the $X_i's$. A standard metric on $G$ , called the Carnot-Caratheodory metric, is naturally associated with $\mathbb X$ and is defined as follows: let $\ell : [0,1] \to G$ be an absolutely continuous path. We say that $\ell$ is admissible if there exist measurable functions $a_1,...,a_k : [0,1] \to \mathbb C$ such that, for almost every $t \in [0,1]$, one has 
$$\ell'(t)=\sum_{i=1}^k a_i(t) X_i(\ell(t)).$$
If $\ell$  is admissible, its length is defined by 
$$|\ell |= \int_0^1\left(\sum_{i=1}^k |a_i(t)|^2 \,dt \right)^{ \frac 12 }.$$

For all $x,y \in G $, define $d(x,y)$ as  the infimum of the lengths  of all admissible paths joining $x$ to $y$ (such a curve exists by the H\"ormander condition). This distance is left-invariant. For short, we denote by $|x|$ the distance between $e$, the neutral element of the group and $x$,  so that the distance from $x$ to $y$ is equal to  $|y^{-1}x|$. 

For all $r>0$, denote by $B(x,r)$ the open ball in $G$ with respect to the Carnot-Caratheodory distance and by $V(r)$ the Haar measure of any ball. 
%There exists $d\in \N^{\ast}$ (called the local dimension of $(G\mathbb X)$) and $0<c<C$ such that, for all $r\in (0,1)$,
%$$
%cr^d\leq V(r)\leq Cr^d,
%$$
%see \cite{nsw}. When $r>1$, two situations may occur (see \cite{guivarch}): 
%\begin{itemize}
%\item Either there exist $c,C,D >0$ such that, for all $r>1$, 
%$$c r^D \leq V(r) \leq C r^D$$
%where $D$ is called the dimension at infinity of the group (note that, contrary to $d$, $D$ does not depend on $\mathbb X$). The group is said to have polynomial volume growth. 
%\item Or  there exist $c_1,c_2,C_1,C_2 >0$ such that, for all $r>1$, 
%$$c_1 e^{c_2r} \leq V(r) \leq C_1 e^{C_2r}$$
%and the group is said to have exponential volume growth. 
%\end{itemize} 
%When $G$ has polynomial volume growth, it is plain to see that there exists $C>0$ such that, for all $r>0$,
%\begin{equation} \label{homog}
%V(2r)\leq CV(r)
%\end{equation}
%which implies that there exist $C>0$ and $\kappa>0$ such that, for all $r>0$ and all $\theta>1$,
%\begin{equation} \label{homogiter}
%V(\theta r)\leq C\theta^{\kappa}V(r).
%\end{equation}
%\bigskip

We denote 
$$-\Delta_G  =-\sum_{i=1}^k X_i^2$$
the sub-laplacian on $G$ and $\nabla_G=(X_1,...,X_k)$ the associated gradient.

%We are interested in proving the existence of precised Hardy inequalities on $G$. When $G=\mathbb H^d$ is the Heisenberg group, such inequalities were proved in \cite{BCG} using para-product. More precisely, it is proved that for all $s \in (0,N/2]$, there exists a constant $C$ such that 
%$$\int_{\mathbb H^d} \frac{|u|^2}{|x|^{2s}}\,dx \leq C \|u\|^2_{ \dot{H}^s(\mathbb H^d)},$$
%where $ \dot{H}^s(\mathbb H^d)$ is defined via complex interpolation. 

%
%The purpose of this paper is two-fold: simplify with the proof of the previous inequality and consider the case of general lie groups. More precisely, the proof of the previous inequality relies on para-product arguments which were developped in \cite{BG}. In the present paper, we provide  a direct proof without the use of para-product techniques. Furthermore, we provide such an inequality in the case of $H-$type groups, which are particular Carnot groups and finally provide a general result for lie groups with polynomial growth. This latter result is very general and relies on the existence of a standard Hardy inequality.

The Lie group $G$ is  called a Carnot group if and only if $G$ is simply connected and the Lie algebra of $G$ admits a stratification, i.e. there exist linear subspaces $V_1,...,V_k$ of $\mathcal G$ such that 
$$\mathcal G= V_1 \oplus ... \oplus V_k$$
which
$$[V_1,V_i]=V_{i+1}$$
for $i=1,...,k-1$ and $[V_1,V_k]=0$. By $[V_1,V_i]$ we mean the subspace of $\mathcal G$ generated by the elements $[X,Y]$ where $X \in V_1$ and $Y \in V_i$. Recall that the class of Carnot groups is a strict subclass of nilpotent groups. Moreover, if $G$ is a Carnot group, there exists $n\in \N$, called the homogeneous dimension of $G$, such that, for all $r>0$,
\begin{equation} \label{vol}
V(r)\sim r^n
\end{equation}
(see \cite{FS}).  The Heisenberg group ${\mathcal H}^d$ is a Carnot group and $n=2d+2$. \par

\medskip

\noindent Let $G$ be a Carnot group, denote by $\delta$ the Dirac distribution supported at the origin and let $u$ be a solution of 
$$-\Delta_G u=\delta.$$
Define $N(x)=u(x)^{\frac{1}{2-n}}$ for $x\neq 0$ and $N(0)=0$. The function $N$ is an homogeneous norm on $N$ by \cite{fol}. Kombe \cite{kombe} proved the following Hardy inequality on $G$: let $\alpha>2-n$ and $u\in {\mathcal D}(G\setminus \left\{0\right\})$, then the following holds 
\begin{equation}\label{standardCarnot}
\left(\frac{n+\alpha-2}2\right)^2\int_G u^2(x)\frac{|\nabla_G N(x)\vert^2}{|N(x)|^2}N^{\alpha}(x)dx\leq  \int_G  \left\vert \nabla_G u(x)\right\vert^2 N^{\alpha}(x)dx.
\end{equation}
We prove the following non local version of \eqref{standardCarnot}:  
\begin{theo} \label{mainth1} 
Let $G$ be a Carnot group with homogeneneous dimension $n \geq 3$. Then for all $\alpha>2-n$ and all $s \in (0,1)$,
 \begin{equation} \label{hardycarnotfrac}
\int_{G} u^2(x)
    \left( \frac{|\nabla_G N\vert}{|N|}\right)^{s}N^{\alpha}dx \lesssim  \iint_{G\times G} \frac{\left\vert u(x)-u(y)\right\vert^2}{\left\vert
       y^{-1} x\right\vert^{n+s}}N^{\alpha}(x) \, dx\, dy \ \forall u\in {\mathcal D}(G\setminus \left\{0\right\}). 
\end{equation}
\end{theo}
\subsection{The case of Riemannian manifolds}
A general principle was developed in \cite{carron} to derive Hardy inequalities on Riemannian manifolds. Let us recall here an example of such an inequality. Let $M$ be a Riemannian manifold, denote by $n$ its dimension, by $d\mu$ its Riemannian measure, by $d$ the exterior differentiation and by $\Delta$ the Laplace-Beltrami operator. For all $x\in M$ and all $r>0$, let $B(x,r)$ be the open geodesic ball centered at $x$ with radius $r$, and $V(x,r)$ its measure. \par
\noindent Assume that $\rho:M\rightarrow [0,+\infty)$ satisfies
\begin{equation} \label{drho}
\left\vert d\rho\right\vert\leq 1\mbox{ on }M,
\end{equation}
and
\begin{equation} \label{subsol}
\Delta\rho\leq -\frac C{\rho}\mbox{ in the distribution sense},
\end{equation}
where $C>0$. Then, for all $\alpha<C-1$ and all $u\in {\mathcal D}(M\setminus \rho^{-1}(0))$,
\begin{equation} \label{hardyriem}
\left(\frac{C-1-\alpha}2\right)^2\int_M \left(\frac u{\rho}\right)^2(x)\rho^{\alpha}(x)dx\leq \int_M \left\vert du(x)\right\vert^2\rho^{\alpha}(x)dx.
\end{equation}
Moreover, if the codimension of $\rho^{-1}(0)$ is greater than $2-\alpha$, \eqref{hardyriem} holds for all function $u\in {\mathcal D}(M)$ (see \cite{carron}, Th\'eor\`eme 1.4 and Remarque 1.5) . \par

\medskip

\noindent We provide here a non local version of \eqref{hardyriem}. To state this result, we introduce some extra assumptions on $M$. The first one is a Faber-Krahn inequality on $M$. For any bounded open subset $\Omega\subset M$, denote by $\lambda_1^D(\Omega)$ the principal eigenvalue of $-\Delta$ on $\Omega$ under the Dirichlet boundary condition. If $p\geq n$, consider the following Faber-Krahn inequality: there exists $C>0$ such that
\begin{equation} \label{FK}
\lambda_1^D(\Omega)\geq C\mu(\Omega)^{\frac 2p}\mbox{ for all bounded subset }\Omega\subset M.
\end{equation}
Let $\Lambda_p>0$ be the greatest constant for which \eqref{FK} is satisfied. In other words,
$$
\Lambda_p=\inf \frac{\lambda_1^D(\Omega)}{\mu(\Omega)^{\frac 2p}},
$$
where the infimum is taken over all bounded subsets $\Omega\subset M$. The Faber-Krahn inequality \eqref{FK} is satisfied in particular when an isoperimetric inequality holds on $M$: namely there exists $C>0$ and $p\geq n$ such that, for all bounded smooth subset $\Omega\subset M$,
\begin{equation} \label{isop}
\sigma(\partial\Omega)\geq C\mu(\Omega)^{1-\frac 1p},
\end{equation}
where $\sigma(\partial\Omega)$ denotes the surface measure of $\partial\Omega$. If $M$ has nonnegative Ricci curvature, then \eqref{isop} with $p=n$ and \eqref{FK} with $p=n$ are equivalent. More generally, if $M$ has Ricci curvature bounded from below by a constant, \eqref{FK} with $p>2n$ implies \eqref{isop} with $\frac p2$ (\cite{carronsmf}, Proposition 3.1, see also \cite{coulhon} when the injectivity radius of $M$ is furthermore assumed to be bounded). Note that there exists a Riemannian manifold satisfying \eqref{FK} for some $p\geq n$ but for which \eqref{isop} does not hold for any $p\geq n$ (\cite{carronsmf}, Proposition 3.4).  \par

\medskip

\noindent It is a well-known fact that \eqref{FK} implies a lower bound for the volume of geodesic balls in $M$. Namely (\cite{carronsmf}, Proposition 2.4), if \eqref{FK} holds, then, for all $x\in M$ and all $r>0$,
\begin{equation} \label{minvol}
V(x,r)\geq \left(\frac{\Lambda_p}{2^{p+2}}\right)^{\frac p2}r^p.
\end{equation}

\medskip

We will also need another assumption on the volume growth of the balls in $M$. Say that $M$ has the doubling property if and only if there exists $C>0$ such that, for all $x\in M$ and all $r>0$,
\begin{equation}Ê\label{doub} \tag{$D$}
V(x,2r)\leq CV(x,r).
\end{equation}
There is a wide class of manifolds on which (\ref{doub}) holds. First, 
it is true on Lie groups with polynomial volume growth (in particular 
on nilpotent Lie groups). Next, (\ref{doub}) is true if $M$ has 
nonnegative Ricci curvature thanks to the Bishop comparison theorem 
(see \cite{BC}). Recall also that (\ref{doub}) remains valid if $M$
is quasi-isometric to a manifold with nonnegative Ricci curvature, or
is a cocompact covering manifold  whose deck transformation group has
polynomial
growth, \cite{CSC}. Contrary to the doubling property, the 
nonnegativity of the Ricci curvature is not stable under 
quasi-isometry. \par

\medskip

We prove the following theorem. 
\begin{theo} \label{hardyriemnonloc}
Let $M$ be a complete non compact Riemannian manifold. Assume that \eqref{FK} holds and that $M$ has the doubling property. Assume also that $C>0$ and $\rho:M\rightarrow [0,+\infty)$ are such that \eqref{drho} and \eqref{subsol} hold. Then, if $\alpha<C-1$ and $\rho^{-1}(0)$ has codimension greater than $2-\alpha$, one has, for all $s\in (0,1)$,\begin{equation} \label{hardyineqriem}
\int_M u^2(x)\rho^{s(\alpha-2)}(x)dx\lesssim \iint_{M\times M} \frac{\left\vert u(y)-u(x)\right\vert^2}{d(x,y)^{p+s}} \rho^{\alpha}(x)dxdy\ \forall u\in {\mathcal D}(M\setminus \rho^{-1}(0)).
\end{equation}
\end{theo}

%\begin{rem} \label{rem1} 
%Analogous ideas show that, if $G$ is  a  unimodular lie group with polynomial growth and if the following inequality holds 
%\begin{eqnarray} \label{standardHardy}
%\int _{G  } |\nabla f |^2
%    \geq  \lambda 
%   \, \int_{G}\frac{ \left\vert
%     f(x)\right\vert^2 }{|x|^{2}}dx,
%\end{eqnarray}

%the following self-improvement holds: if $s\in
%(0,2)$, there exists $\lambda>0$  such that, for any function $f\in {\mathcal D}(G)$

% \begin{eqnarray} \label{poincfrac} \iint _{G \times G }
%   \frac{\left\vert f(x)-f(y)\right\vert^2}{V\left(\left\vert
%       y^{-1} x\right\vert\right)\left\vert
%       y^{-1} x\right\vert^{s}} \, dx\, dy \ge
%   \lambda '
%   \, \int_{G}\frac{ \left\vert
%     f(x)\right\vert^2 }{|x|^{2s}}dx. 
%\end{eqnarray}
%\end{rem}

\section{Proof of Theorem \ref{mainth1}}

In order to prove Theorem \ref{mainth1}, we need to introduce an operator $L_{N^{\alpha}}$ on $L^2(G)$. Let $L^2(G,N^{\alpha})$ denote the $L^2$ space on $G$ with respect to the measure $N^{\alpha}dx$ and $H^1(G,N^{\alpha})$ the Sobolev space defined by
$$
H^1(G,N^{\alpha}):=\left\{f\in L^2(G,N^{\alpha});\ X_if\in L^2(G,N^{\alpha})\ \forall 1\leq i\leq k\right\}.
$$
Define now the operator $L_{N^{\alpha}}$ on $L^2(G,N^{\alpha})$ by
$$
L_{N^{\alpha}}u:=-N^{-\alpha}\sum_{i=1}^k X_i(N^{\alpha}X_iu).
$$
The domain of $L_{N^{\alpha}}$ is given by
$$
{\mathcal D}\left(L_{N^{\alpha}}\right)=\left\{u\in L^2(G,N^{\alpha});\ N^{-\alpha} X_i(N^{\alpha}X_iu)\in L^2(G,N^{\alpha})\ \forall 1\leq i\leq k\right\}.
$$
One has, for all $u\in {\mathcal D}\left(L_{N^{\alpha}}\right)$ and all $v\in H^1(G,N^{\alpha})$,
$$
\int_G L_{N^{\alpha}}u(x)v(x)N^{\alpha}(x)dx=\int_G \sum_{i=1}^k X_iu(x)X_iv(x)N^{\alpha}(x)dx.
$$
The operator $L_{N^{\alpha}}$ is therefore clearly symmetric and nonnegative on $L^2(G,N^{\alpha})$, and the spectral theorem allows to define the
usual powers $(L_{N^{\alpha}})^{\beta}$ for any $\beta>0$ by means of spectral
theory. 

By the definition of $L_{N^{\alpha}}$,  (\ref{standardCarnot}) means, in terms of operators  in $L^2(G,N^{\alpha})$, that, for some $\lambda>0$,
\begin{equation} \label{ineqop}
L_{N^{\alpha}} \geq \lambda \mu, 
\end{equation}
where $\mu$ is the multiplication operator by $\frac{|\nabla_G N\vert}{|N|}$. Using a functional calculus argument (see \cite{davies}, p.
110) one deduces from (\ref{ineqop}) that, for any $s\in (0,2)$,
\[
(L_{N^{\alpha}})^{s/2}\geq \lambda^{s/2}\mu^{s/2}
\]
which implies, thanks to the fact $(L_{N^{\alpha}})^{s/2} = ((L_{N^{\alpha}})^{s/4})^2$ and the
symmetry of $(L_{N^{\alpha}})^{s/4}$ on $L^2(G,N^{\alpha})$, that
$$
\int_{G} \left\vert u(x)\right\vert^2\mu(x)^sN^{\alpha}(x)dx \leq  C \int_{G} \left\vert (L_{N^{\alpha}})^{s/4}u(x)\right\vert^2N^{\alpha}(x)dx =$$
$$
C \left\| (L_{N^{\alpha}})^{s/4} u \right\|^2 _{L^2(G,N^{\alpha})}.
$$
The conclusion of Theorem \ref{mainth1} follows now from the estimate
$$ 
\left\Vert (L_{N^{\alpha}})^{s/4} u \right\|^2 _{L^2(G)} \leq C \, \iint_{G \times G}
\frac{\left\vert u(x)-u(y)\right\vert^2}{\left\vert y^{-1}x\right\vert^{n+2s}} N^\alpha(x)    \, dx\, dy,
$$
which is exactly the estimate for $M=N^{\alpha}$ provided in Lemmata 3.2 and 3.3 in \cite{rslie} (remember that \eqref{vol} holds).

\section{Proof of Theorem \ref{hardyriemnonloc}}
The proof of Theorem \ref{hardyriemnonloc} relies again on estimates for the powers of a suitable operator. Namely, define
$$
L_{\rho^{\alpha}}u=\rho^{-\alpha}\mbox{div}(\rho^{\alpha}\nabla u),
$$
where $\nabla$ is the gradient induced by the Riemannian metric and $\mbox{div}$ is the divergence operator on $M$. As before, $L_{\rho^{\alpha}}$ is a nonnegative symmetric operator on $L^2(M,\rho^{\alpha}dx)$ and $L_{\rho^{\alpha}}^{\beta}$ is defined by spectral theory for all $\beta>0$. If $\mu$ denotes the multiplication operator by $\rho^{-2}$, \eqref{hardyriem} means that $L_{\rho^{\alpha}}\geq c\mu$ in $L^2(M,\rho^{\alpha}dx)$. Spectral theory then yields $L_{\rho^{\alpha}}^{s/2}\geq c\mu^{s/2}$, which means that
$$
\int_M u^2(x)\rho^{s(\alpha-2)}(x)dx\lesssim \left\Vert L_{\rho^{\alpha}}^{s/4}u\right\Vert_{L^2(M,\rho^{\alpha}dx)}^2,
$$
and we are therefore left with the task of checking
$$
\left\Vert L_{\rho^{\alpha}}^{s/4}u\right\Vert_{L^2(M,\rho^{\alpha}dx)}^2\lesssim \iint_{M\times M} \frac{\left\vert u(y)-u(x)\right\vert^2}{d(x,y)^{p+s}} \rho^{\alpha}(x)dxdy.
$$
To that purpose, one first notices that
$$
\left\Vert L_{\rho^{\alpha}}^{s/4}u\right\Vert_{L^2(M,\rho^{\alpha}dx)}^2\lesssim \int_0^{+\infty} t^{-1-s/2} \left\Vert t \, L_{\rho^{\alpha}} \,
(\mbox{I} + t \, L_{\rho^{\alpha}})^{-1} u\right\Vert_{L^2(M,\rho^{\alpha}dx)}^2 \, dt
$$
and it is therefore enough to show that
$$
\int_0^{+\infty} t^{-1-s/2} \left\Vert t \, L_{\rho^{\alpha}} \,
(\mbox{I} + t \, L_{\rho^{\alpha}})^{-1} u\right\Vert_{L^2(M,\rho^{\alpha}dx)}^2 \, dt $$
$$\lesssim \iint_{M\times M} \frac{\left\vert u(y)-u(x)\right\vert^2}{d(x,y)^{p+s}} \rho^{\alpha}(x)dxdy.
$$
The proof follows the same lines as the one of  Lemma 3.3 in \cite{rslie} and we will therefore be sketchy, only indicating the main differences. Using (\ref{doub}), one can pick up a countable family $x_j^t$, $j\in \N$, such that the balls
$B\left(x_j^{t},\sqrt{t}\right)$ are pairwise disjoint and 
\begin{equation} \label{union}
M=\bigcup_{j \in \N} B\left(x_j^{t},2\sqrt{t}\right).
\end{equation}
By \eqref{doub}, there exist constants $\tilde C>0$ and $\kappa>0$
such that for all $\theta>1$ and all $x\in G$, there are at most $\tilde
C\, \theta^{2\kappa}$ indexes $j$ such that $|x^{-1}x_j^t |\leq
\theta\sqrt{t}$.

\medskip

For fixed $j$, one has
$$
t \, L_{\rho^{\alpha}} \, (\mbox{I} + t\, L_{\rho^{\alpha}})^{-1} u= t \, L_{\rho^{\alpha}} \, (\mbox{I} + t\, L_{\rho^{\alpha}})^{-1} \,
g^{j,t}
$$
where, for all $x\in M$, 
$$
g^{j,t}(x):=u(x)-m^{j,t}
$$
and $m^{j,t}$ is defined by
$$
m^{j,t}:=\frac 1{V\left(x_j^t,2\sqrt{t}\right)}\int_{B\left(x_j^{t},2\sqrt{t}\right)}
u(y) dy.$$
Note that, here, the mean value of $u$ is computed with respect to the
Riemannian measure on $M$. Since (\ref{union}) holds, one
clearly has
$$
\begin{array}{lll}
\displaystyle \left\Vert t \, L _{\rho^{\alpha}}\, (\mbox{I} + t\, L_{\rho^{\alpha}})^{-1} u 
\right\Vert_{L^2(M,\rho^{\alpha}dx)}^2 
& \leq & \displaystyle \sum_{j \in \N} 
\left\Vert t \, L_{\rho^{\alpha}} \, (\mbox{I} + t\, L_{\rho^{\alpha}})^{-1} u
\right\Vert_{L^2\left(B(x_j^t,2\sqrt{t}),\rho^{\alpha}dx\right)}^2\\
& = & \displaystyle  \sum_{j \in \N} \left\Vert t\, L_{\rho^{\alpha}} \, (\mbox{I} + t\,
L_{\rho^{\alpha}})^{-1}  g^{j,t}
\right\Vert_{L^2\left(B\left(x_j^t,2\sqrt{t}\right),\rho^{\alpha}dx \right )}^2,
\end{array}
$$
and it is therefore enough to ckeck
\begin{eqnarray}\label{laststep}
\sum_{j\in \N} \left\Vert t \, L_{\rho^{\alpha}} \, (\mbox{I} + t\, L_{\rho^{\alpha}})^{-1}
g^{j,t}\right\Vert_{L^2\left(B\left(x_j^{t},2\sqrt{t}\right),\rho^{\alpha}dx \right )} ^2\\\nonumber
\lesssim  \iint_{M\times M} \frac{\left\vert u(y)-u(x)\right\vert^2}{d(x,y)^{p+s}} \rho^{\alpha}(x)dxdy.
\end{eqnarray}
As in \cite{rslie}, this is a consequence of $L^2$ off-diagonal estimates for $L_{\rho^{\alpha}}$ and upper estimates for the functions $g^{j,t}$. Let us recall the off-diagonal estimates for $L_{\rho^{\alpha}}$ for completeness:
\begin{lem} \label{off} There exists $C$ with the following property: for all closed disjoint subsets
$E,F\subset M$ with $\mbox{d}(E,F)=:d>0$, all function $f\in
L^2(M,\rho^{\alpha}dx)$ supported in $E$ and all $t>0$,
$$
\left\Vert (\mbox{I}+t \, L_{\rho^{\alpha}})^{-1}f\right\Vert_{L^2(F,\rho^{\alpha}dx)}+\left\Vert t \,
 L_{\rho^{\alpha}}(\mbox{I}+t \, L_{\rho^{\alpha}})^{-1}f\right\Vert_{L^2(F,\rho^{\alpha}dx)}\leq $$
 $$ 8 \, e^{-C \,
 \frac{d}{\sqrt t}} \left\Vert f\right\Vert_{L^2(E,\rho^{\alpha}dx)}.
$$
\end{lem}
The proof of Lemma \ref{off} is analogous to the one of Lemma 3.1 in \cite{rslie}. \par

\medskip

\noindent As far as estimates for $g^{j,t}$ are  concerned, set, for all $k \ge 1,$
$$
C_0^{j,t}=B\left(x_j^t, 4\sqrt{t}\right) \ \mbox{ and
} \ C_k^{j,t}=B\left(x_j^t,2^{k+2}\sqrt{t}\right)\setminus
B\left(x_j^t,2^{k+1}\sqrt{t}\right), 
$$
and $g^{j,t}_k:=g^{j,t} \, {\bf 1}_{C_k^{j,t}}$, $k \ge 0$, where, for any
subset $A\subset M$, ${\bf 1}_A$ is the usual characteristic function of
$A$. We then have:
\begin{lem} \label{estimg}
There exists $\bar C>0$ such that, for all $t>0$ and all $j \in \N$:
\begin{itemize}
\item[{\bf A.}] 
$$\displaystyle \left\Vert g_0^{j,t}\right\Vert_{L^2(C_0^{j,t},\rho^{\alpha}dx)}^2\leq  
\frac{\bar C}{t^{p/2}} \int_{B\left(x_j^t,4\sqrt{t}\right)}
\int_{B\left(x_j^t,4\sqrt{t}\right)} \left\vert u(x)-u(y)\right\vert^2 \,
\rho^{\alpha}(x)dx \, dy. $$
\item[{\bf B.}]
For all $k\geq 1$,
\[ \left\Vert g^{j,t}_k\right\Vert_{L^2(C_k^{j,t},\rho^{\alpha}dx)}^2 
\leq \]
\[
\frac{\bar C}{(2^k\sqrt{t})^p} \int_{x\in B(x^t_j,2^{k+2}\sqrt{t})} 
\int_{y\in B(x^t_j,2^{k+2}\sqrt{t})} \left\vert u(x)-u(y)\right\vert^2 \, \rho^{\alpha}(x)dx\, dy.\]
\end{itemize}
\end{lem}
The proof of Lemma \ref{estimg} is analogous to the one of Lemma 3.4 in \cite{rslie}, the only extra ingredient being the lower bound \eqref{minvol} applied with the balls $B(x^t_j,2^{k+2}\sqrt{t})$.  We then conclude the proof of \eqref{laststep} in the same way as for the conclusion of the proof of Lemma 3.3 in \cite{rslie}. 
\bibliographystyle{alpha}                  % *.bst files
   \bibliography{Biblio-hardy-riem}      

   \medskip  
{\em Emmanuel Russ}--
Universit\'e Paul C\'ezanne, LATP,\\
Facult\'e des Sciences et Techniques, Case cour A\\
Avenue Escadrille Normandie-Niemen, F-13397 Marseille, Cedex 20, France et  \\
CNRS, LATP, CMI, 39 rue F. Joliot-Curie, F-13453 Marseille Cedex 13, France 

\medskip 

{\em Yannick Sire}--
Universit\'e Paul C\'ezanne, LATP,\\
Facult\'e des Sciences et Techniques, Case cour A\\
Avenue Escadrille Normandie-Niemen, F-13397 Marseille, Cedex 20, France et  \\
CNRS, LATP, CMI, 39 rue F. Joliot-Curie, F-13453 Marseille Cedex 13, France.

\end{document}